\documentclass{article}
\usepackage{amsmath,amssymb,,amsfonts,rsfso,pictex}

\def\cl#1{{\mathscr #1}}
\newcommand{\dlines}{\displaylines}
\def\limn{\lim_{n\to\infty}}
\def\ep{\varepsilon}

\def\span{\mathop{\rm span}}
\def\eqref#1{(\ref{#1})}
\def\proof{\noindent{\it Proof}}

\newcommand{\R}{\mathbb{R}}
\renewcommand{\th}{\theta}

\newcommand{\e}{{\rm e}}

\def\tin#1{\par\noindent\hskip3em\llap{#1\enspace}\ignorespaces}
\def\cl#1{{\mathcal #1}}

\def\E{{\rm E}}

\def\P{{\rm P}}

\def\FINEDIM{\par\hfill$\blacksquare$\hphantom{mm}
\medskip

\noindent }
\def\tfrac#1#2{{\textstyle\frac {#1}{#2}}}
\def\supp{\mathop{\rm supp}}

\newtheorem{theorem}{Theorem}[section]

\newtheorem{rem}{Remark}[section]

\newtheorem{lem}{Lemma}[section]

\begin{document}
\title{\bf Intermediate spaces, Gaussian probabilities and exponential tightness}
\author{Paolo Baldi\protect\footnote{baldi@mat.uniroma2.it}\protect\footnote{The author wishes to thank G. Pisier for valuable suggestions and advice.}\protect\footnote{The author wishes to thank one of the referees whose remarks have given rise to significant improvement of the paper.}\protect\footnote{The author acknowledges the MIUR Excellence Department Project awarded to the Dipartimento di Matematica, Universit\`a  di  Roma  ``Tor  Vergata'',  CUP  E83C18000100006}\\
Dipartimento di Matematica, Universit\`a di Roma "Tor Vergata"\\
}
\medskip


\date{}\maketitle
\begin{abstract}
We prove the existence of an intermediate Banach space between the space where the Gaussian measure lives and its RKHS, thus extending what happens with Wiener measure, where the intermediate space can be chosen as a space of H\"older paths. From this result it is very simple to deduce a result of exponential tightness for Gaussian probabilities. 
\end{abstract} 
\bigskip

\noindent{\it AMS 2000 subject classification:} 60F10, 60B12
\smallskip

\noindent{\it Key words and phrases:} Gaussian probabilities, interpolation theory.

\section{Introduction}
Let $E=\cl C_0([0,T],\R^m)$ be the space of continuous $\R^m$-valued paths starting at $0$ and endowed with the sup norm and let $\mu$ be the Wiener measure on it. It is well known that the Reproducing Kernel Hilbert Space (RKHS) of this Gaussian probability is the space $\cl H=H^1_0([0,T])$ of the paths $\gamma$ vanishing at $0$, that are absolutely continuous and have a square integrable derivative. 

Let us denote by $\cl C^0_\alpha\subset E$ the space of the paths that are $\alpha$-H\"older continuous and whose modulus of continuity 
$$
\omega(\delta)=\sup_{{0\le s<t\le T\atop|t-s|\le\delta}}|\gamma(t)-\gamma(s)|
$$
is such that
$$
\lim_{\delta\to 0+}\frac {\omega(\delta)}{\delta^\alpha}=0\ .
$$
This space (whose elements are sometimes called the ``small'' $\alpha$-H\"older continuous functions) is separable and it is also well known (thanks to Kolmogorov's continuity theorem) that, for $0<\alpha<\frac 12$, $\mu(\cl C^0_\alpha)=1$ (whereas $\mu(\cl H)=0$). It is also well known that, still for $0<\alpha<\frac 12$,
$$
E\hookleftarrow \cl C^0_\alpha\hookleftarrow\cl H
$$
{\it the embeddings being compact}.

We are concerned with the question whether the existence of such an ``intermediate'' space is a general fact, i.e. true for every centered Gaussian probability on a separable Banach space.

More precisely we prove the following result.
\begin{theorem}\label{main} Let $E$ be a separable Banach space, $\mu$ a centered Gaussian probability on $E$ and $\cl H$ the corresponding RKHS. Then there exists a Banach space $\widetilde E$,  separable and such that 
\tin{a)} $\mu(\widetilde E)=1$ and
\tin{b)} the embeddings 
$$
E\hookleftarrow\widetilde E\hookleftarrow\cl H
$$
are compact.
\end{theorem}
We shall even prove that there are infinitely many such spaces. We shall call ``intermediate space'' any separable Banach space satisfying  a) and b) of Theorem \ref{main}.

The proof of Theorem \ref{main} is the object of \S\ref{sec-main}. Of course Theorem \ref{main} is obvious if $E$ is finite dimensional, as then we can choose $E=\widetilde E=\cl H$. Therefore in the sequel we implicitly assume that $E$ is infinite dimensional.  

In the proof we shall also assume that $E=\supp(\mu)$. Otherwise just consider $\supp (\mu)$ instead of $E$.

This investigation was motivated by  an application to the Large Deviations of the sequence of probabilities $(\mu_\ep)_\ep$, $\mu_\ep$ being the image of $\mu$ through the map $x\mapsto\ep x$. The property of exponential tightness is a key step in the proof of these estimates. One remarks that its proof in the case of Wiener measure is particularly simple and is based, besides Fernique's theorem, on the existence of the spaces of H\"older continuous functions, which are intermediate spaces. Thanks to Theorem \ref{main} the same, simple proof of exponential tightness for the Wiener measure works for a general Gaussian probability on a separable Banach space as developed in \S\ref{sec-lg}. 

\S\ref{sec-comm} is devoted to comments and complements.

The proof of Theorem \ref{main} is largely inspired to the arguments of the fundamental papers of L.Gross \cite{gross-TAMS} and \cite{gross-berkeley}.
\section{Motivation: exponential tightness of Gaussian measures}\label{sec-lg}
Let $\mu$ be a centered Gaussian probability on the separable Banach space $E$ and let $\mu_\ep$ be its image through the map $x\mapsto \ep x$ as above. 

The Large Deviations properties of the family $\mu_\ep$ as $\ep\to0$  are well understood since a long time (see see \cite{dvIII} \S 5 e.g.). One of the main steps in this investigation is to prove that the family $(\mu_\ep)_\ep$ is {\it exponentially tight} at speed $\ep\mapsto\ep^2$, i.e. that for every $R>0$ there exists a compact set $K_R\subset E$ such that
\begin{equation}\label{eq-et}
\limsup_{\ep\to 0}\ep^2\log \mu_\ep(K_R^c)\le -R\ .
\end{equation}
This fact follows immediately from Theorem \ref{main}: let us denote $\|\enspace\|_i$ the norm of the intermediate space $\widetilde E$, as $\|\enspace\|_i$ is $\mu$-a.s. finite (this is a) of Theorem \ref{main}), by Fernique's theorem (\cite{fernique-cras}, \cite{fernique-stflour74}), for some $\rho>0$, we have
\begin{equation*}
\int_E\e^{\rho\|x\|_i^2}\,d\mu(x):=C_\rho<+\infty\ .
\end{equation*}
Let $K_R$ denote the ball of radius $\sqrt{R/\rho}$ of $\widetilde E$, which is compact in $E$. 
From 
$$
C_\rho\ge \int_{K_R^c/\ep}\e^{\rho\|x\|_i^2}\,d\mu(x)\ge 
\mu(\tfrac 1{\ep}K_R^c)\e^{R/\ep^2}
$$
we deduce
$$
\limsup_{\ep\to 0}\ep^2\log \mu_\ep(K_R^c)=\limsup_{\ep\to 0}\ep^2\log \mu(\tfrac 1{\ep} K_R^c)\le-R\ ,
$$
i.e. \eqref{eq-et}.
\section{Proof of the main result}\label{sec-main}
From now on $\mu$ will denote a Gaussian probability on the infinite dimensional separable Banach space $E$ as in the introduction, $\|\enspace\|$ being the norm of the Banach space $E$.

If we denote by $E'$ the dual of $E$, then to every continuous functional $\xi\in E'$ we can associate the r.v. $(E,\mu)\to \R$ defined  as
$x\mapsto\langle \xi,x\rangle$. Let $E'_\mu$ be the completion of $E'$ in $L^2(\mu)$. This is a separable Hilbert space which is also a Gaussian space. 

For every $g\in E'_\mu$ the vectors 
\begin{equation}\label{eq-alter}
h=\int_E xg(x)\, d\mu(x)
\end{equation}
form a vector space $\cl H\subset E$ which, endowed with the scalar product
$$
\langle h_1,h_2\rangle_{\cl H}=\int_E g_1(x)g_2(x)\, d\mu(x)
$$
for 
$$
h_1=\int_E xg_1(x)\, d\mu(x),\qquad h_2=\int_E xg_2(x)\, d\mu(x)\ ,
$$
is an Hilbert space $\cl H$ isometric to $E'_\mu$. Remark that \eqref{eq-alter} can also be written $h=\E[Xg(X)]$, $X$ denoting an $E$-valued r.v. having law equal to $\mu$. $\cl H$ is called the Reproducing Kernel Hilbert Space (RKHS) of $\mu$.

For more details on the structure of Gaussian probabilities  see \cite{ledoux-talagrand} or the very nice and very short presentation in \S 2 of \cite{deacosta-small}.

%
\medskip

\noindent{\it Proof} of Theorem \ref{main}. Recall that we assume $E$ to be infinite dimensional. 

\tin{a)} First step: construction of $\widetilde E$.

Let $(\Omega,\cl F,\P)$ be a probability space and $X:\Omega\to E$ a Gaussian r.v. having distribution $\mu$. Let $(g_n)_n$ be an orthonormal system of the Hilbert space $E'_\mu$. This forms also a sequence of independent $N(0,1)$-distributed r.v.'s. Let $e_n=\E[X g_n(X)]$. $e_n\in E$ and $(e_n)_n$ is an orthonormal system of the RKHS $\cl H$.

Then it is well known (see Proposition 3.6 p. 64 in \cite{ledoux-talagrand}) that the sequence
\begin{equation}\label{eq-Xn}
X_n=\sum_{j=1}^ng_je_j 
\end{equation}
is a square integrable $E$-valued martingale converging a.s. and in $L^2$ to $X$. Hence
$$
\limn \E\Bigl(\Bigl\|\sum_{j=n}^\infty g_je_j\Bigr\|^2 \Bigr)=0\ .
$$
Let $\alpha>0$ be fixed and let $(n_k)_k$ be an increasing sequence of integers such that $n_0=0$ and
\begin{equation}\label{eq-nk}
\E\Bigl(\Bigl\|\sum_{j=n_k+1}^{\infty} g_je_j\Bigr\|^2 \Bigr)\le 2^{-k(3+2\alpha)}\ .
\end{equation}
Let $\cl H_k=\span(e_{n_k+1},\dots,e_{n_{k+1}})$ and let $Q_k$ be the projector $\cl H\to\cl H_k$. Let, for the vector $x=\sum_{n=1}^\infty \alpha_n e_n\in \cl H$,
\begin{equation}\label{eq-rinf}
\|x\|_i:=\sum_{k=0}^\infty 2^{k\alpha}\Big\|\sum_{j=n_k+1}^{n_{k+1}}\alpha_j e_j\Big\|=
\sum_{k=1}^\infty 2^{k\alpha}\|Q_k(x)\|\ .
\end{equation}
$\|\enspace\|_i$ is a norm on $\cl H$. Actually Lemma \ref{lem-finitenorm} below states that $\|x\|_i<+\infty$ for every $x\in\cl H$ while subadditivity and positive homogeneity are immediate.
Let 
\begin{equation}\label{eq-wk}
W_k:=\sum_{j=n_k+1}^{n_{k+1}}g_j e_j\ .
\end{equation}
The $E$-valued r.v.'s $W_k$ are Gaussian and independent. Remark that, with our choice of the numbers $n_k$, thanks to the Markov inequality we have
\begin{equation}\label{key-ineq}
\P(2^{k\alpha}\|W_k\|\ge 2^{-k})\le 2^{2k(1+\alpha)}\,\E(\|W_k\|^2)\le2^{2k(1+\alpha)}\,\E\Bigl(\Bigl\|\sum_{j=n_k+1}^{\infty} g_je_j\Bigr\|^2 \Bigr)\le 2^{-k}\ .
\end{equation}
We can now define $\widetilde E=$the completion of $\cl H$ with respect to the norm $\|\enspace\|_i$. Remark that, as for $x\in\cl H$
$$
\|x\|_i=
\sum_{k=1}^\infty 2^{k\alpha}\|Q_k(x)\|\ge 
\sum_{k=1}^\infty \|Q_k(x)\|\ge\Big\|\sum_{k=1}^\infty Q_k(x)\Big\| =\|x\|\ ,
$$
we have $\widetilde E\subset E$. It is also obvious that $\widetilde E$ is dense in $E$, as it contains $\cl H$ which is itself dense in $E$.
\tin{b)} Second step: $\mu(\widetilde E)=1$. Let 
$$
Y_k=X_{n_k}=\sum_{j=1}^{n_k}g_j e_j
$$
and let us prove that $(Y_k)_k$, as a sequence of $\widetilde E$-valued r.v.'s, converges in probability. We proceed quite similarly as in the proof of the subsequent Lemma \ref{lem-finitenorm}. Let $\ep>0$. If $p_0$ is such that $2^{-p_0}<\frac \ep2$, then for $p_0\le \ell< r$,
$$
\dlines{
\P(\|Y_r-Y_\ell\|_i>\ep)\le\P\Bigl(\sum_{k=\ell+1}^{r}2^{\alpha k}\|W_k\|>\sum_{k=\ell+1}^r2^{-k}\Bigr)\le\cr \le\sum_{k=\ell+1}^r\P\bigl(2^{\alpha k}\|W_k\|>2^{-k}\bigr)\le \sum_{k=p_0+1}^{\infty}2^{- k}\le2\cdot 2^{-p_0}=\ep\ .\cr
}
$$
Therefore $(Y_k)_k$ is a Cauchy sequence in probability in $\widetilde E$, hence it converges, in probability, to some $\widetilde E$-valued r.v. $\widetilde X$. As $\widetilde E\subset E$ and its topology is stronger, $(Y_k)_k$ also converges in probability in $E$. But we know already that $(Y_k)_k$ in $E$ converges to a r.v. $X$ having law $\mu$. Hence $\widetilde X=X$ a.s. and 
$\mu(\widetilde E)=\P(X\in \widetilde E)=\P(\widetilde X\in \widetilde E)=1$.

\tin{c)} Step three: the embedding $E\hookleftarrow\widetilde E$ is compact. 

Let $(x_p)_p$ be a bounded sequence in $\widetilde E$ and $(z_p)_p\subset \cl H$ another sequence such that $\|x_p-z_p\|_i<2^{-p}$, which is possible as $\cl H$ is dense in $\widetilde E$. Let $M$ be such that $\|z_p\|_i\le M$ for every $p$. As the projectors $Q_k$ are finite dimensional, for every $k$ there exists a subsequence $(p_r^{(k)})_r$ such that $\|Q_kz_{p_r^{(k)}}-y^{(k)}\|_i\to 0$ for some vector $y^{(k)}\in \cl H_k$, i.e. of the form
$$
y^{(k)}=\sum_{m=n_k+1}^{n_{k+1}} \alpha_m e_m\ .
$$
By the diagonal argument there exists a subsequence $(p'_r)_r$ such that 
$\|Q_kz_{p'_r}-y^{(k)}\|_i\to 0$ as $r\to\infty$ for every $k$. Let now
$\ep>0$ be fixed. We have, for every positive integer $k_0$,
\begin{equation}\label{eq-diagonal}
\|z_{p'_r}-z_{p'_\ell}\|=\Bigl\|\sum_{k=1}^\infty Q_k(z_{p'_r}-z_{p'_\ell})\Bigr\|\le \Bigl\|\sum_{k=1}^{k_0} Q_k(z_{p'_r}-z_{p'_\ell})\Bigr\|+\Bigl\|\sum_{k=k_0+1}^\infty Q_k(z_{p'_r}-z_{p'_\ell})\Bigr\|\ .
\end{equation}
We first choose $k_0$ so that $M\,2^{-\alpha k_0}<\frac \ep3$, so that
$$
\dlines{
\Bigl\|\sum_{k=k_0+1}^\infty Q_k(z_{p'_r}-z_{p'_\ell})\Bigr\|\le
\sum_{k=k_0+1}^\infty \Bigl\|Q_k(z_{p'_r}-z_{p'_\ell})\Bigr\|\le\cr
\le2^{-\alpha k_0}\sum_{k=k_0+1}^\infty 2^{\alpha k} \Bigl\|Q_k(z_{p'_r}-z_{p'_\ell})\Bigr\|\le 2^{-\alpha k_0}\|z_{p'_r}-z_{p'_\ell}\|_i\le 2^{-\alpha k_0}\cdot 2M\le\frac23\, \ep\cr
}
$$
and then $p_0$ so that, for $r,\ell\ge p_0$
$$
\Bigl\|\sum_{k=1}^{k_0} Q_k(z_{p'_r}-z_{p'_\ell})\Bigr\|\le \frac \ep3\ \cdotp
$$
Therefore $(z_{p'_r})_r$ is a Cauchy sequence in $E$, hence also $(x_{p'_r})_r$ which proves the compactness of the embedding $\widetilde E\hookrightarrow E$.
\tin{d)} Last step $\widetilde E \hookleftarrow \cl H$ is compact. This is immediate as, $\widetilde E$ being dense in $E$, $\cl H$ is also the RKHS of the Gaussian probability $\mu$ on $\widetilde E$. Such an embedding is always compact.
\FINEDIM
\begin{lem}\label{lem-concentration} Let $\cl K$ be a finite dimensional Hilbert space and $\cl F\subset\cl K$ a subspace and let us denote by $\nu$ and $\nu'$ the standard $N(0,I)$ distributions on $\cl K$ and $\cl F$ respectively. Let $B\subset \cl K$ be a convex, centrally symmetric, convex set. Then
$$
\nu(B)\le\nu'(\cl F\cap B)\ .
$$
\end{lem}
For the proof of Lemma \ref{lem-concentration} the reader is directed to \cite{gross-TAMS} as this is a weaker version of Lemma 4.1 there.
\begin{lem}\label{lem-finitenorm}
$\|x\|_i<+\infty$ for every $x\in\cl H$.
\end{lem}
\proof Let us first prove that the sequence of real r.v.'s 
$$
Z_n=\sum_{k=1}^n 2^{k\alpha}\|W_k\|
$$
converges in probability. 
Let $\ep>0$ and $p_0$ such that $2^{-p_0}<\frac \ep2$. Remark that this implies $\ep>\sum_{k=p_0+1}^\infty2^{-k}$. For $p_0\le n<m$, we have thanks to \eqref{key-ineq}
$$
\dlines{
\P(|Z_n-Z_m|>\ep)=\P\Bigl(\sum_{k=n+1}^m 2^{k\alpha}\|W_k\|>\ep\Bigr)\le
\P\Bigl(\sum_{k=n+1}^m 2^{k\alpha}\|W_k\|>\sum_{k=n+1}^m2^{-k}\Bigr)\le \cr
\le\sum_{k=n+1}^m\P\Bigl( 2^{k\alpha}\|W_k\|>2^{-k}\Bigr)\le 
\sum_{k=n+1}^m 2^{-k}\le \sum_{k=p_0+1}^\infty 2^{-k}=2\cdot 2^{-p_0}< \ep\ .\cr
}
$$
Hence $(Z_n)_n$ is a Cauchy sequence in probability 
and converges in probability to some real r.v. $Z$. 

Let us prove that, for every $\ep>0$, we have $\P(Z<\ep)>0$. Let
$$
Z_N=\sum_{k=1}^{N}2^{k\alpha}\|W_k\|,\qquad Z-Z_N=\sum_{k=N+1}^\infty 2^{k\alpha}\|W_k\|\ .
$$
We have $\P(Z_N<\frac\ep2)>0$, as $Z_N$ depends only on the modulus of finitely many r.v.'s $W_k$ each of them being Gaussian and with values in a finite dimensional vector space. Moreover let $N$ be large enough so that
$\P(Z-Z_N<\frac\ep2)>0$. As the r.v.'s $Z-Z_N$ and $Z_N$ are independent (they depend on different $W_k$'s)
\begin{equation}\label{eq-ep}
\P(Z<\ep)\ge \P\Bigl(Z-Z_N<\frac\ep2,Z_N<\frac\ep2\Bigr)=\P\Bigl(Z-Z_N<\frac\ep2\Bigr)
\P\Bigl(Z_N<\frac\ep2\Bigr)>0\ .
\end{equation}
Finally, let us assume that there exists $x=\sum_{i=1}^\infty \alpha_je_j\in \cl H$ such that $\|x\|_i=+\infty$ and let us prove that this is absurd. Of course we can assume $\|x\|_{\cl H}=1$. Let us consider the seminorms on $\cl H$
$$
\|z\|_k=\sum_{j=1}^k 2^{j\alpha}\|Q_j(z)\|
$$
so that $\lim_{k\to\infty} \|z\|_k=\|z\|_i$. Let $\cl K_k=\span(e_1,\dots,e_{n_k})$ so that the r.v. $X_{n_k}$ of \eqref{eq-Xn} takes values in $\cl K_k$. Let $x^{(k)}=\sum_{i=1}^{n_k} \alpha_je_j$ be the projection of $x$ on $\cl K_k$. Remark that $\|x\|_k=\|x^{(k)}\|_k$. We apply Lemma \ref{lem-concentration} considering the convex set 
$$
B=\{z\in\cl K_k, \|z\|_k\le a\}
$$
and $\cl F=\span(x^{(k)})$. Let $\xi_k=\sqrt{\alpha_1^2+\dots+\alpha_{n_k}^2}$, so that the vector 
$
\frac {x^{(k)}}{\xi_k}
$
has modulus $1$ in $\cl H$ and the r.v.
$$
g:=\frac 1{\xi_k}\sum_{j=1}^{n_k}\alpha_jg_j
$$
is $N(0,1)$-distributed. Hence the r.v.
$$
\frac1{\xi_k}\sum_{j=1}^{n_k} \alpha_jg_j\cdot\frac {x^{(k)}}{\xi_k}
$$
is $N(0,1)$-distributed and $\cl F$-valued.

Let $a$ be a continuity point of the partition function of $Z$. Thanks to Lemma \ref{lem-concentration}, as $\xi_k\to \|x\|_{\cl H}=1$ and we assume $\|x\|_i=+\infty$, we have
$$
\dlines{
\P(Z\le a)=\lim_{k\to\infty}\P\Bigl(\Bigl\|\sum_{j=1}^{n_k}g_j e_j\Bigr\|_i\le a\Bigr)\le
\lim_{k\to\infty}\P\Bigl(\Bigl|\frac 1{\xi_k}\sum_{j=1}^{n_k}\alpha_jg_j\Bigl|\frac{\|x\|_k}{\xi_k}\le a\Bigr)=\cr
=\lim_{k\to\infty}\P\Bigl(|g|\le\frac {a\xi_k}{\|x\|_k} \Bigr)=0\cr
}
$$
which is in contradiction with \eqref{eq-ep} and completes the proof of Lemma \ref{lem-finitenorm}.
\FINEDIM
\section{Remarks and complements}\label{sec-comm}%
\begin{rem}\label{rem1}\rm In fact we have proved the existence on infinitely many intermediate spaces $\widetilde E$ between $E$ and $\cl H$ (recall that we assume that $E$ is infinite dimensional). Actually the argument above can be repeated in order to construct a subsequent intermediate space $\widetilde E_1$ between $\cl H$ and $\widetilde E$, which will be necessarily different of $\widetilde E$, the embedding
$\widetilde E\hookleftarrow\widetilde E_1$ being compact. And so on.
\end{rem}

\begin{rem}\label{rem2}\rm In a first attempt to prove Theorem \ref{main} the author tried considering interpolation spaces. More precisely let, for $x\in E$,
$$
K(t,x)=\inf_{a+b=x, a\in E, b\in \cl H}(\|a\|+t|b|_{\cl H})
$$
and let, for $0<\th<1$,
\begin{equation}\label{xf}
\|x\|_\th=\sup_{t>0}t^{-\th}K(t,x)\ .
\end{equation}
Let us define the vector space $G_\th$ as the set of vectors $x\in E$ such that $\|x\|_\th<+\infty$, endowed with the norm $\|\enspace\|_\th$. See \cite{bergh-lofstrom}, \cite{lunardi} or \cite{pisier} for more details on this topic.

It is well known that $G_\th$ is a Banach space and also that the embeddings $E\hookleftarrow G_\th\hookleftarrow\cl H$ are compact (\cite{lions-peetre}, \S V.2). 

The question remains whether $\mu(G_\th)=1$.

In the case of the Wiener space, $E=\cl C([0,T],\R)$, $\cl H=H^1_0$ and $\mu=$the Wiener measure, it can be proved that 
$G_\th$ contains the space of small $\alpha$-H\"older paths for $\alpha>\th$, which is a separable Banach space having Wiener measure $1$. This gives $\mu(G_\th)=1$ for $\th<\frac 12$, which however leaves open the question in the case $\th\ge\frac 12$. 

The author does not know whether such a 
Banach space $G_\th$ is also separable, but it can be proved that the closure of $H^1_0$ in $G_\th$, $\widetilde G_\th$ say, also contains the small $\alpha$-H\"older paths for $\alpha>\th$. Hence $\widetilde G_\th$, which is separable, is an intermediate space in the sense of Theorem \ref{main} in this case. Note that the requirement $\th<\frac 12$ means that $G_\th$ should be ``closer'' to $E$ than to $\cl H$. 

Concerning interpolation spaces, hence, many questions, possibly of interest, remain open.
Is it true, in general, that the interpolated space $G_\th$ is an intermediate space? For every $0<\th<1$ or just for some values of the interpolating parameter $\th$?
\end{rem}
\begin{rem}\label{rem-ref1}\rm The construction of the intermediate space $\widetilde E$ of \S\ref{sec-main} is of course not unique, as other possibilities are available for the candidate norm \eqref{eq-rinf}. For instance, let us define the sequence $(n_k)_k$ so that
\begin{equation}\label{eq-nk2}
\E\Bigl(\Bigl\|\sum_{j=n_k+1}^\infty g_je_j\Bigr\|^2 \Bigr)\le 2^{-2k(\alpha+\eta)}\ .
\end{equation}
for some $\eta>0$ and then, for $x=\sum_{n=1}^\infty \alpha_n e_n\in \cl H$,
\begin{equation}\label{eq-rinf2}
\|x\|'=\sup_{k>0}2^{k\alpha}\Big\|\sum_{j=n_k+1}^{n_{k+1}}\alpha_j e_j\Big\|\ .
\end{equation}
In this remark we prove that 
also $\|x\|'<+\infty$ for $x\in\cl H$ and that the completion of $\cl H$ with respect to $\|\enspace\|'$ is also an intermediate space.

This has some interest as it shows that there are (many) other possible ways of constructing intermediate spaces. Also we shall see, in the next remark, that for a suitable choice of the orthonormal system $(e_n)_n$, in the case $E=\cl C([0,1])$ and Wiener measure the resulting intermediate spaces are the H\"older spaces.

The proof of $\|x\|'<+\infty$ for $x\in\cl H$ is actually even simpler than the one of Lemma \ref{lem-finitenorm}: let $W_k$ as in \eqref{eq-wk} and for $n>0$
$$
Z_n=\sup_{k\le n} 2^{k\alpha}\|W_k\|\ .
$$
Let us show that $Z_n\to Z$ where $Z$ is a r.v. such that $\P(Z<\ep)>0$ for every $\ep>0$. The a.s. convergence of $(Z_n)_n$ is immediate being an increasing sequence. Denoting by $Z$ its limit we have
$$
\dlines{
\P(Z\le \ep)=\limn\P\bigl(2^{k\alpha}\|W_k\|<\ep, k=1,\dots,n\bigr)=\cr
=\limn\prod_{k=1}^n\P(2^{k\alpha}\|W_k\|\le \ep)=\prod_{k=1}^\infty\bigr(1-\P(2^{k\alpha}\|W_k\|>\ep)\bigl)\ .\cr
}
$$
The infinite product above converges to a strictly positive number if and only if the series $\sum_{k=1}^\infty\P(\|W_k\|>\ep)$ is convergent. But, by Markov's inequality and using the bound \eqref{eq-nk2},
\begin{equation}\label{eq-convprime1}
\P(2^{k\alpha}\|W_k\|>\ep)\le \frac 1{\ep^2}\,2^{2k\alpha}\E(\|W_k\|^2)\le \frac 1{\ep^2}\,2^{-2k\eta}
\end{equation}
which is the general term of a convergent series. Moreover the limit $Z$ is finite a.s., as the r.v.'s $W_k$ are independent and the event $\{Z=+\infty\}$ is a tail event having probability $<1$.

The remainder of the proof of Lemma \ref{lem-finitenorm} is quite similar to the one developed in \S\ref{sec-main}. 

The new norm $\|\enspace\|'$ also produces a family of intermediate spaces, as stated in the next result
\begin{theorem}\label{main2} Let $\widetilde E'=$the completion of $\cl H$ with the norm $\|\enspace\|'$. Then $\widetilde E'$ is an intermediate space.
\end{theorem}
\proof Let us prove first that $\mu(\widetilde E')=1$. Let, as in the proof of Theorem \ref{main},
$Y_k=X_{n_k}=\sum_{j=1}^{n_k}g_j e_j$
 and let us prove that $(Y_k)_k$, as a sequence of $\widetilde E'$-valued r.v.'s, converges in probability. If $k_0\le\ell\le r$ we have now
$$
\P(\|Y_r-Y_\ell\|_i>\ep)=\P\Bigl(\sup_{\ell\le k\le r}2^{\alpha k}\|W_k\|>\ep\Bigr)\le\P\Bigl(\sup_{k\ge k_0}2^{\alpha k}\|W_k\|>\ep\Bigr)\ .
$$
Markov's inequality gives $\P(2^{\alpha k}\|W_k\|>\ep)\le\frac 1{\ep^2}\, 2^{-k\eta}$, hence by the Borel-Cantelli Lemma $2^{\alpha k}\|W_k\|>\ep$ for finitely many $k$ only and
$$
\lim_{k_0\to\infty}\P\Bigl(\sup_{k\ge k_0}2^{\alpha k}\|W_k\|>\ep \Bigr)=0
$$
so that $(Y_k)_k$ is a Cauchy sequence in probability in $\widetilde E'$ which implies, with the same argument as in the proof of Theorem \ref{main}, that $\mu(\widetilde E')=1$.

We are left with the proof that the embedding $\widetilde E'\hookleftarrow\widetilde E$ is compact, which is quite similar to the argument of the proof of Theorem \ref{main}. 

Let again $(x_p)_p$ be a bounded sequence in $\widetilde E$ and $(z_p)_p\subset \cl H$ another sequence such that $\|x_p-z_p\|'<2^{-p}$, which is possible as $\cl H$ is dense in $\widetilde E$. Let $M$ be such that $\|z_p\|'\le M$ for every $p$. As the projectors $Q_k$ are finite dimensional, for every $k$ there exists a subsequence $(p_r^{(k)})_r$ such that $\|Q_kz_{p_r^{(k)}}-y^{(k)}\|'\to 0$ for some vector $y^{(k)}\in \cl H_k$, i.e. of the form
$y^{(k)}=\sum_{m=n_k+1}^{n_{k+1}} \alpha_m e_m$.
By the diagonal argument there exists a subsequence $(p'_r)_r$ such that 
$\|Q_kz_{p'_r}-y^{(k)}\|_i\to 0$ as $r\to\infty$ for every $k$. Let now
$\ep>0$ be fixed. We have for every positive integer $k_0$
$$
\|z_{p'_r}-z_{p'_\ell}\|=\Bigl\|\sum_{k=1}^\infty Q_k(z_{p'_r}-z_{p'_\ell})\Bigr\|\le \Bigl\|\sum_{k=1}^{k_0} Q_k(z_{p'_r}-z_{p'_\ell})\Bigr\|+\Bigl\|\sum_{k=k_0+1}^\infty Q_k(z_{p'_r}-z_{p'_\ell})\Bigr\|\ .
$$
As $\|Q_k(z_{p'_r}-z_{p'_\ell}))\|\le2\cdot2^{-k\alpha}\|Q_k(z_{p'_r}-z_{p'_\ell}))\|'\le 2\cdot2^{-k\alpha} M$, we have
$$
\Bigl\|\sum_{k=k_0+1}^\infty Q_k(z_{p'_r}-z_{p'_\ell})\Bigr\|\le \sum_{k=k_0+1}^\infty \bigl\|Q_k(z_{p'_r}-z_{p'_\ell})\bigr\|\le
2 M\sum_{k=k_0+1}^\infty2^{-k\alpha}\le  2^{-k_0\alpha}\frac {2 M}{1-2^{-\alpha}}
$$
which, for some $k_0$ large enough is $\le \frac\ep2$. We can choose now $p_0$ so that, for $r,\ell\ge p_0$
$$
\Bigl\|\sum_{k=1}^{k_0} Q_k(z_{p'_r}-z_{p'_\ell})\Bigr\|\le \frac \ep2\ \cdotp
$$
Therefore $(z_{p'_r})_r$ is a Cauchy sequence in $E$, hence also $(x_{p'_r})_r$ which proves the compactness of the embedding $\widetilde E\hookrightarrow E$.
\FINEDIM

Note that the norms $\|\enspace\|_i$ and $\|\enspace\|'$ differ not only because of their definitions (\eqref{eq-rinf} as opposed to \eqref{eq-rinf2}) but also on the different requirements on the sequence $(n_k)_k$ (\eqref{eq-nk} and \eqref{eq-nk2}).
\end{rem}
One of the referees raised the natural question whether the intermediate spaces $\widetilde E$ constructed in Theorem \ref{main}, in the case of the Wiener space might produce the H\"older spaces, or, more generally, if they can be described in terms of regularity of the paths, with the idea that the larger the parameter $\alpha$, the greater the regularity of the paths of $\widetilde E$.

The question in general seems to the author to require an analysis going beyond the scope of the present paper, in particular taking into account that regularity also depends on the choice of the orthonormal system $(e_n)_n$ and possibly on the regularity of its elements. 

In the next remark we show however that, for a certain choice of the orthonormal system, the intermediate spaces of Theorem \ref{main2} are actually the H\"older spaces $\cl C^0_\alpha$.
\begin{rem}\label{rem-ref}\rm Let $E=\cl C([0,1],\R)$ and $\|\enspace\|$ the sup norm.

Let us recall the characterization, due to Ciesielski \cite{Cies1}, of the small Holder spaces $\cl C_\alpha^0$.

Let $\{\chi_n\}_n$ be the Haar system, namely the set of functions on the interval $[0,1]$ defined as $\chi_1(t) \equiv 1$ and
$$
\chi_{2^k+j}(t) =\begin{cases}
\sqrt{2^k}\hfil\qquad &{\rm if\ }t\in [{2j-2\over 2^{k+1}},{2j-1\over
2^{k+1}}[\\
-\sqrt{2^k}\hfil\qquad &{\rm if\ }t\in [{2j-1\over 2^{k+1}},{2j\over
2^{k+1}}[\hfill\cr
0 &{\rm otherwise}\cr
\end{cases}
$$
for $k=1,2,\dots$, $j=1,2,\dots,2^k$. It is well known that $\{\chi_n\}_n$ is a
complete orthonormal system of $L^2([0,1],\R)$. 
\begin{figure}[h!]
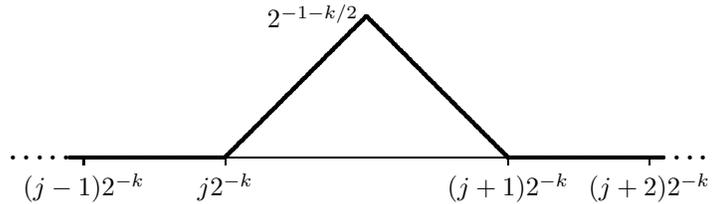

\hbox to \hsize\bgroup\hss
\beginpicture
\setcoordinatesystem units <.37truein,.37truein>
\setplotarea x from -4.2 to 4.2, y from 0 to 2
\axis bottom shiftedto y=0 ticks short  withvalues $(j-1)2^{-k}$ $j2^{-k}$ $(j+1)2^{-k}$ $(j+2)2^{-k}$ / at -4 -2 2 4 / /
\axis left  shiftedto x=0 invisible ticks length <0pt> withvalues $2^{-1-k/2}$ / at  2 / /
\setplotsymbol ({\bf.})
\plot -4.2 0 -2 0  0 2  2 0 4.2 0 /
\setdots
\plot -5 0 -4.2 0 /
\plot 4.2 0  5 0 /
\endpicture
\hss\egroup
\caption{The graph of $\phi_{2^k+j}$.\label{ex3.sommadiuniform}}
\end{figure}
Moreover let $\phi_n(t)=\int_0^t\chi_n(s)\, ds$ be the
primitive of $\chi_n$ (the Schauder basis).  For a continuous path $x\in \cl C([0,1],\R)$ let us consider the coefficients  $\xi_n=\int_0^1\chi_n(s)\,dx(s)$
which are well defined, as $\chi_n$ is piecewise constant. 
Ciesielski \cite{Cies1} proved that the separable Banach spaces $\cl C_\alpha^0$ and $c_0$ (the sequences vanishing at $\infty$ endowed with the sup norm) are isomorphic ($0<\alpha< 1$). 

More precisely if
$$
w_{2^k+j}(\alpha)=2^{k(\alpha-\frac 12)+(1-\alpha)}
$$
then $x\in\cl C_\alpha^0,\ 0<\alpha< 1$, if and only if
$\xi=\{\xi_nw_n(\alpha)\}_n\in c_0$. Let us denote by $c_\alpha$ the space of the sequences $\xi=(\xi_n)_n$ such that
$(\xi_nw_n(\alpha))_n\in c_0$.
Ciesielski's theorem states that the mapping
$$
(\xi_n)_n\enspace \leftrightarrow \sum_{m=1}^\infty \xi_m\phi_m
$$
is an isomorphism between $c_\alpha$ and $\cl C^0_\alpha$.

Remark that under Ciesielski's isomorphism $\cl H$ is mapped into $\ell_2$. Actually  $(\phi_n)_n$ is itself an orthonormal  basis of $\cl H=H^1_0$.

The spaces $\cl C^0_\alpha$ for $0<\alpha<\frac 12$ are actually the intermediate spaces obtained in Remark \ref{rem-ref1}, if we choose, as an orthonormal system for $\cl H$, $e_n=\phi_n$. Actually, if $n_k=2^k$ we have, noting that, the supports of the $\phi_j$ are disjoint and 
$\|\phi_j\|= 2^{-1-k/2}$ for $2^k+1\le j\le 2^{k+1}$, for large $k$
$$
\E\Bigl(\Bigl\|\sum_{j=2^k+1}^{2^{k+1}} g_j\phi_j\Bigr\|^2 \Bigr)=2^{-2-k}\E\Bigl(\sup_{2^k+1\le j\le 2^{k+1}}g_j^2\Bigr)\le 2^{-k\lambda}
$$
for every $\lambda<1$. This is an elementary, but a bit involved, computation that we are not going to explicit here. Hence, as in Remark \ref{rem-ref1}, the norm
$$
\|x\|'=\sup_{k>0}2^{k\alpha}\Big\|\sum_{j=n_k+1}^{n_{k+1}}\alpha_j e_j\Big\|
$$
is finite on $\cl H$ as soon as $2(\alpha+\eta)<1$, i.e. $\alpha<\frac 12$. 
Note, again using the fact that the supports of the $\phi_j$ are disjoint and 
$\|\phi_j\|= 2^{-1-k/2}$ for $2^k+1\le j\le 2^{k+1}$, 
$$
\dlines{
\sup_{k>0}2^{k\alpha}\Big\|\sum_{j=2^k+1}^{2^{k+1}}\xi_j\phi_j\Big\|=
\sup_{k>0}\sup_{2^k+1\le j\le 2^{k+1}}2^{k\alpha}2^{-1-k/2} |\xi_j| = const\,|\xi|_\alpha \cr
}
$$
so that, thanks to Ciesielski's isomorphism, the intermediate norm $\|\enspace\|'$ is finite if and only if $x\in \cl C^0_\alpha$.
\end{rem}

\def\cprime{$'$} \def\cprime{$'$} \def\cprime{$'$}
\providecommand{\bysame}{\leavevmode\hbox to3em{\hrulefill}\thinspace}
\providecommand{\MR}{\relax\ifhmode\unskip\space\fi MR }
\providecommand{\MRhref}[2]{%
  \href{http://www.ams.org/mathscinet-getitem?mr=#1}{#2}
}
\providecommand{\href}[2]{#2}

\end{document}